\documentclass[12pt,reqno]{amsart}

\usepackage{verbatim}
\setbox0=\hbox{$+$}
\newdimen\plusheight
\plusheight=\ht0
\def\+{\;\lower\plusheight\hbox{$+$}\;}

\setbox0=\hbox{$-$}
\newdimen\minusheight
\minusheight=\ht0
\def\-{\;\lower\minusheight\hbox{$-$}\;}

\setbox0=\hbox{$\cdots$}
\newdimen\cdotsheight
\cdotsheight=\plusheight%\ht0
\def\cds{\lower\cdotsheight\hbox{$\cdots$}}

\makeatletter
\def\leqalignno#1{\displ@y \tabskip\z@ plus\@ne fil
  \halign to\displaywidth{\hfil$\@lign\displaystyle{##}$\tabskip\z@skip
    &$\@lign\displaystyle{{}##}$\hfil\tabskip\z@ plus\@ne fil
    &\kern-\displaywidth\rlap{$\@lign\hbox{\rm##}$}\tabskip\displaywidth\crcr
    #1\crcr}}
\makeatother

\newcommand{\ra}{\rightarrow}
\newcommand{\eb}{\begin{equation}}
\newcommand{\ee}{\end{equation}}

\newcommand{\df}{\dfrac}
\newcommand{\tf}{\tfrac}

\newcommand{\Q}{\mathbb{Q}}
\newcommand{\C}{\mathbb{C}}
\newcommand{\N}{\mathbb{N}}

\renewcommand{\(}{\left\(}
\renewcommand{\)}{\right\)}
\renewcommand{\[}{\left\[}
\renewcommand{\]}{\right\]}
\renewcommand{\i}{\infty}
\renewcommand{\pmod}[1]{\,(\textup{mod}\,#1)}
\numberwithin{equation}{section}
 \theoremstyle{plain}
\newtheorem{theorem}{Theorem}[section]

\newtheorem{corollary}[theorem]{Corollary}

\begin{document}

\title[Determinations of Analogues of Gauss Sums]
{Determinations of Analogues of Gauss Sums and Other Trigonometric
Sums}
\author{Matthias Beck, Bruce C.~Berndt1, O-Yeat Chan, and Alexandru
  Zaharescu}
\address{Department of Mathematics, San Francisco State University, 1600
Holloway Ave, San Francisco, CA 94132}
\address{Department of Mathematics, University of Illinois, 1409 West Green
Street, Urbana, IL 61801, USA}
\email{beck@math.sfsu.edu}\email{berndt@math.uiuc.edu}
\email{ochan@math.uiuc.edu} \email{zaharesc@math.uiuc.edu}

%\maketitle

\vspace*{0.5in}
\begin{center}

{\bf DETERMINATIONS OF ANALOGUES OF GAUSS SUMS AND OTHER
TRIGONOMETRIC
  SUMS}\\[5mm]

{\footnotesize MATTHIAS BECK, BRUCE
C.~BERNDT\footnotemark[1]\footnotemark[2], O-YEAT
CHAN\footnotemark[3], and
ALEXANDRU ZAHARESCU}\\[3mm]

\end{center}

\vspace*{0.5in}

\footnotetext[1]{Research partially supported by  grant
MDA904-00-1-0015   from the National Security Agency.}

\footnotetext[2]{Research partially supported by the Korea
Institute for Advanced Study.}

\footnotetext[3]{Research partially supported by a grant from the
University of Illinois Research Board.}

\vskip 5mm \noindent{\footnotesize{\bf Abstract.} Explicit
determinations of several classes of trigonometric sums are given.
These sums can be viewed as analogues or generalizations of Gauss
sums.  In a previous paper, two of the present authors  considered
primarily sine sums associated with primitive odd characters. In
this paper, we establish two general theorems involving both sines
and cosines, with more attention given to cosine sums in the
several examples that we provide.}

\vskip 2mm \noindent{\footnotesize {\bf 2000 AMS Classification
Numbers}:
  Primary, 11L03;   Secondary, 11R29, 11L10.}

\section{Introduction}

Motivated by two trigonometric identities
\begin{equation}\label{ident1}
\df{\sin(2\pi/7)}{\sin^2(3\pi/7)}
-\df{\sin(\pi/7)}{\sin^2(2\pi/7)}+
\df{\sin(3\pi/7)}{\sin^2(\pi/7)} = 2\sqrt{7}
\end{equation}
and
\begin{equation}\label{ident2}
\df{\sin^2(3\pi/7)}{\sin(2\pi/7)}
-\df{\sin^2(2\pi/7)}{\sin(\pi/7)}+
\df{\sin^2(\pi/7)}{\sin(3\pi/7)} = 0,
\end{equation}
discovered by Berndt  and L.--C.~Zhang \cite{bz} as corollaries of
two theta-function identities in Ramanujan's notebooks \cite{nb},
and also motivated by further identities found by Z.--G.~Liu
\cite[pp.~107--108]{zhiguo}, Berndt and Zaharescu \cite{bz1}
evaluated large classes of trigonometric sums in terms of class
numbers of imaginary quadratic fields.  As an illustration, we
begin by offering one of these general theorems and one of its
corollaries.

\begin{theorem}\label{t1}
Let $\chi$ denote an odd, real, nonprincipal, primitive character
of modulus $k$, where $k$ is odd and $k \geq 7$. Let
\begin{equation}\label{2.1}
S_1(k,\chi,a,b) := \sum_{0<n<k/2}\chi(n)\df{\sin^a(b\pi
n/k)}{\sin^{a+1}(\pi n/k)},
%S_1(k) := \sum_{0<n<k/2}\chi(n)\df{\sin(4\pi n/k)}{\sin^2(\pi
%n/k)}.
\end{equation}
where $a$ and $b$ are positive integers with $a$ odd and $b$ even.
%and $ab-a-3 < 2k$.
Define
 \begin{equation}\label{C}
  C_{a,b} := \sum_{\substack{n,m,r \geq 0 \\
n+bm+r=(ab-a-1)/2}\\}(-1)^m\chi(n)\binom{a}{m}\binom{a+r}{r}.
\end{equation}
 Then
\begin{equation}\label{2.2}
S_1(k,\chi,a,b) =\sqrt{k}\left(b^a\,h(-k)-2C_{a,b}\right),
\end{equation}
where $h(-k)$ denotes the class number of the imaginary quadratic
field $\mathbb{Q}(\sqrt{-k})$.
\end{theorem}

% The requirement $ab-a-3 <2k$ is technical; we could
%dispense with this condition, but then $C_{a,b}$ would have to be
%replaced by a more complicated sum.

\begin{corollary}\label{c0} If $\chi$ is given as above, then
\begin{equation}\label{c1}
\sum_{0<n<k/2}\chi(n)\df{\sin(4\pi n/k)}{\sin^{2}(\pi n/k)} =
\sqrt{k}\left(4h(-k)-2\right).
\end{equation}
\end{corollary}

The evaluation \eqref{ident1} is the special case of \eqref{c1}
when $k=7$ and $\chi(n)$ is the Legendre symbol
$\left(\tfrac{n}{7}\right)$.

The proofs in \cite{bz1} depend on contour integration and
elementary properties of Gauss sums.

While sums involving powers of the sine function, odd characters,
and class numbers of imaginary quadratic fields are the focus of
\cite{bz1}, in this paper we consider sums with both sines and
cosines, and sums involving either odd or even characters. Some of
our main theorems are for odd characters, and others are for even
characters.
%Although some of our evaluations can be expressed
%in terms of class numbers in special cases, we do not emphasize
%such instances in the present paper.
Class numbers arise in our results involving odd characters.

Before embarking on the proofs of our general theorems, we thought
it best to begin with the evaluation of a class of cosine sums
associated with even characters.  Thus, in Section \ref{s2}, we
alter our ideas from \cite{bz1} to evaluate a large class of
trigonometric sums involving even characters.

In Section \ref{main}, we establish our main general theorem,
while in the following section we apply the aforementioned results
 and prove two principal general theorems on sums of
trigonometric functions.  We also offer several corollaries.
%Some of the identities found by Liu \cite{zhiguo} are not associated with characters.
In Section \ref{s3}, we evaluate a large class of trigonometric
sums which includes four of Liu's \cite{zhiguo} identities.

We close the introduction by recording those properties of Gauss
sums that are used in the sequel. Throughout this paper, $\chi$
denotes a nonprincipal, real, primitive character modulo $k$,
where $k$ is an odd positive integer at least equal to 3. Define
the Gauss sum $G(z,\chi)$ for any complex number $z$ by
\begin{equation}\label{gauss}
G(z,\chi) := \sum_{j=0}^{k-1}\chi(j)e^{2\pi ijz/k}.
\end{equation}
Then, for each integer $n$, we have the factorization theorem
\cite[p.~9, Thm.~1.1.3]{bew}
\begin{equation}\label{factor}
G(n,\chi) = \chi(n)G(1,\chi) =:\chi(n)G(\chi).
\end{equation}
In fact, \eqref{factor} characterizes real primitive characters,
i.e., \eqref{factor} holds if and only if $\chi$ is real and
primitive \cite[p.~482, Thm.~1]{apostol},
\cite[pp.~65--66]{daven}.  We need Gauss's famous evaluation
\cite[p.~349, Thm.~7]{bs}, \cite[p.~22, Thm.~1.3.4]{bew}
\begin{equation}\label{gausssum}
G(\chi) =\begin{cases}\sqrt{k}, \qquad &\text{if $\chi$ is even}\\
i\sqrt{k}, \qquad &\text{if $\chi$ is odd}. \end{cases}
\end{equation}

For our results involving class numbers of imaginary quadratic
fields, denote the class number of the field $\Q(\sqrt{-k})$ by
$h(-k)$, and recall the classical formula for the class number
\cite[p.~299]{ayoub}, \cite[p.~344, eq.~(4.3)]{bs}
\begin{equation}\label{class}
h(-k) = -\df{1}{k}\sum_{j=1}^{k-1} j\chi(j),
\end{equation}
which holds when $k\geq 7$ and $\chi$ is odd.

 In the sequel, $R_\alpha(f) = R_\alpha$ denotes the
residue of a meromorphic function $f$ at a pole $\alpha$.

\section{Trigonometric Sums Associated with Even
Characters}\label{s2}

At the outset, we remark that we consider less general functions
in Theorem \ref{thm1}  than we did for the corresponding theorem
in \cite{bz1}, because otherwise the computations of residues
would have been more cumbersome and the results less elegant.

\begin{theorem}\label{thm1}
Let $\chi$ be a real, nonprincipal, even, primitive character of
modulus $k$, where $k$ is odd. For each nonnegative integer $a$
and even positive integer $b$, define
\begin{equation}\label{e1}
S_1(a,b,\chi) := \sum_{0<n<k/2}\chi(n)\df{\cos^a(b\pi
n/k)}{\cos^2(\pi n/k)}
\end{equation}
and
\begin{equation}\label{e2}
E(a,b,\chi) := \df{1}{2^{a-2}}\sum_{\substack{n,j,r \geq 0 \\
2(n+j+br)=ab}}(-1)^jj\chi(n)\binom{a}{r},
\end{equation}
where the sum is over all nonnegative integers $n$, $j$, and $r$
satisfying the condition $2(n+j+br)=ab$. Also, set
\begin{equation}\label{e3}
g(\chi) := \sum_{j=1}^{k-1}(-1)^jj\chi(j).
\end{equation}
Then,
%provided $ab-4 < 2k$,
\begin{equation}\label{e4}
S_1(a,b,\chi) =
-\sqrt{k}\left((-1)^{ab/2}g(\chi)+E(a,b,\chi)\right).
\end{equation}
\end{theorem}

%The condition $ab-4<2k$ is technical; we could eliminate this
%condition, but then the sum $E(a,b,\chi)$ would need to be
%replaced by a more complicated sum.

\begin{proof}
For $N>0$, let $C_N$  denote the positively oriented indented
rectangle with horizontal sides through $\pm iN$ and vertical
sides through $0$ and $k$.  On the left side of $C_N$, there is a
semicircular indentation $I_0$ of radius less than $1$ centered at
$0$ and to its left.  On the right side of $C_N$, the semicircular
indentation comprises the points $I_0+k$. Consider the meromorphic
function
\begin{equation}\label{e5}
f(z)= \df{G(z,\chi)}{G(\chi)}\df{\cos^a(b\pi z/k)}{\cos^2(\pi
z/k)}\df{1}{e^{2\pi iz}-1}.
\end{equation}
 We integrate $f(z)$ over the
contour $C_N$, on the interior of which the function $f(z)$ has
simple poles (at most) at $ z = 1, 2, \dots , k-1$. Also on the
interior of $C_N$, $f(z)$ has a simple pole at $z=k/2$. This pole
is simple because $k$ is odd and because $G(k/2,\chi)=0$.  To see
this, use \eqref{gauss}, replace $j$ by $k-j$, and recall that
$\chi$ is even.  We thus find that $G(k/2,\chi)=-G(k/2,\chi)$,
i.e., $G(k/2,\chi)=0$.  Lastly, since $G(0,\chi)=0$, $f(z)$ does
\emph{not} have a pole at $z=0$.

We first calculate the residues of the poles of $f(z)$ on the
interior of $C_N$. Using \eqref{factor}, we easily find that for
each positive integer $n$,
\begin{equation}\label{e6}
R_n = \df{G(n,\chi)}{G(\chi)}\df{\cos^a(b\pi n/k)}{\cos^2(\pi
n/k)}\df{1}{2\pi i}
 = \df{\chi(n)}{2\pi i}\df{\cos^a(b\pi n/k)}{\cos^2(\pi n/k)}
\end{equation}
and
\begin{equation}\label{e7}
R_{k-n} = \df{G(k-n,\chi)}{G(\chi)}\df{\cos^a(b\pi
(k-n)/k)}{\cos^2(\pi(k-n)/k)}\df{1}{2\pi i}
 = R_n,
 \end{equation}
 because $\chi$ and $b$ are even. Secondly, using \eqref{gausssum},
 we readily find that
 \begin{align}\label{e8}
 R_{k/2} &= \lim_{z\to k/2} \df{1}{G(\chi)}\df{G(z,\chi)}{\cos(\pi z/k)}
 \df{z-k/2}{\cos(\pi z/k)}\df{\cos^a(b\pi z/k)}{e^{2\pi
 iz}-1}\notag\\
 &=
 -\df{(-1)^{ab/2}\sqrt{k}\,i}{\pi}\sum_{j=1}^{k-1}(-1)^jj\chi(j)\notag\\
 &=-\df{(-1)^{ab/2}\sqrt{k}\,i}{\pi}g(\chi),
 \end{align}
 by \eqref{e3}.
 Hence, by the residue theorem, \eqref{e6}--\eqref{e8}, and
 \eqref{e1},
 \begin{equation}\label{e9}
\int_{C_N}f(z)dz = 2(-1)^{ab/2}\sqrt{k}g(\chi) + 2S_1(a,b,\chi).
\end{equation}

Next, we let $N\to\i$ in order to calculate directly the integral
in \eqref{e9}. Firstly, observe from \eqref{e5} that because $b$
is even, $f(z)$ has period $k$, and so the integrals over the
indented vertical sides of $C_N$ cancel. Secondly, examine the
integral over the top horizontal side $C_{NT} $ of $C_N$. Set
$z=x+iN$, $0\leq x \leq k$, and
\begin{equation}\label{19}
\mu := e^{2\pi iz/k} = e^{-2\pi N/k}e^{2\pi ix/k}.
\end{equation}
Then, by \eqref{gauss},
\begin{align}
\df{G(z,\chi)}{e^{2\pi iz}-1} &=
\left(\sum_{j=0}^{k-1}\chi(j)\mu^j\right) \left(-\sum_{n=0}^\infty
\mu^{kn}\right)
= -\sum_{n=0}^\infty \sum_{j=0}^{k-1} \chi(j)\mu^{kn+j} \nonumber \\
\label{20} &=-\sum_{n=0}^\infty \sum_{j=0}^{k-1}
\chi(kn+j)\mu^{kn+j}
=-\sum_{m=0}^\infty \chi(m)\mu^m, \\
\label{21} \cos^a(b\pi z/k) &=
\left(\df{\mu^{b/2}+\mu^{-b/2}}{2}\right)^a =
\df{1}{2^a}\mu^{-ab/2}(1+\mu^b)^a
=\df{1}{2^a}\mu^{-ab/2}\sum_{r=0}^a\binom{a}{r}\mu^{br},\\
\intertext{and} \label{22}
 \cos^{-2}(\pi z/k) &=\df{4}{\mu^{-1}(1+\mu)^2}
 =4\sum_{j=1}^{\i}(-1)^{j-1}j\mu^{j}.
\end{align}
Thus, $f(z)$ has the form \eb\label{12} f(z) = f(x+iN) =
\sum_{n=-ab/2}^0 c_n\mu^n + \sum_{n=1}^\infty c_n\mu^n \ee for
some constants $c_n$, along the top horizontal side $C_{NT}$.
Observe that we can ignore the terms of the form $c_n\mu^n$ with
$n>0$, since their contributions to the integral of $f(z)$ over
$C_{NT}$ tend to 0 as $N$ tends to $\i$ (recall that $|\mu| =
e^{-2\pi N/k}$). Thus, we truncate \eqref{12}, keeping only the
first sum in \eqref{12}, and integrate termwise.  Since, for
$n\neq 0$, \eb\label{13} \int_k^0 \mu^n dx = 0, \ee
 we find that
\begin{equation}\label{14}
\int_{C_{NT}}f(z)dz = -kc_0.
\end{equation}
Since $f(z)$ is an odd function with period $k$, we find that
\begin{align}
\int_{C_{NT}} f(z)dz &= \int_{z=k+iN}^{z=iN} f(z) dz
= -\int_{-z=-k-iN}^{-z=-iN} f(-z)dz \notag\\
%&= \int_{-z=-k-iN}^{-z=-iN} f(-z)d(-z) \notag\\
&= \int_{u=-k-iN}^{u=-iN} f(u)du = \int_{u=-iN}^{u=k-iN} f(u)du
=\int_{C_{NB}} f(z)dz,\label{15}
\end{align}
where $C_{NB}$ is the lower horizontal path of $C_N$. Thus, by
\eqref{e9} and \eqref{14}, we conclude that
\begin{equation}\label{18}
-2kc_0 = 2(-1)^{ab/2}\sqrt{k}g(\chi) + 2S_1(a,b,\chi),
\end{equation}
where $c_0$ is defined in \eqref{12}.  To compute $c_0$, we
utilize  \eqref{20}--\eqref{22} and the definition \eqref{e5} of
$f$ to find that along $C_{NT}$,
\begin{align}\label{24}
f(z)
&=-\df{1}{G(\chi)}\df{\mu^{-ab/2}}{2^{a-2}}\sum_{n=0}^{\infty}\chi(n)\mu^n
\sum_{r=0}^a\binom{a}{r}\mu^{br}\sum_{j=1}^{\i}(-1)^{j-1}j\mu^j.
\end{align}
The constant term in \eqref{24} is equal to
\begin{equation}\label{27}
c_0 = \df{1}{2^{a-2}G(\chi)}
\sum_{\substack{n,j,r \geq 0 \\
n+j+br=ab/2}\\}(-1)^j j\chi(n)\binom{a}{r} =
\df{1}{\sqrt{k}}E(a,b,\chi),
\end{equation}
by \eqref{gausssum} and \eqref{e2}. Hence, from \eqref{18} and
\eqref{27}, we deduce that
\begin{equation*}
-2kc_0 = -2\sqrt{k}E(a,b,\chi)= 2(-1)^{ab/2}\sqrt{k}g(\chi) +
2S_1(a,b,\chi),
\end{equation*}
which completes the proof of Theorem \ref{thm1}.
\end{proof}

We state the special case $a=1$ as a separate corollary.

\begin{corollary}\label{cor1}
 Let $b$ be a positive even integer, and assume
that $\chi$ satisfies the conditions of Theorem
\textup{\ref{thm1}}. Let
\begin{equation}\label{e10}
E(b,\chi) := 2\sum_{\substack{n,j \geq 1\\2n+2j=b}}(-1)^jj\chi(n).
\end{equation}
Then
%if $b<2k+4$,
\begin{equation}\label{e11}
\sum_{0<n<k/2}\chi(n)\df{\cos(b\pi n/k)}{\cos^2(\pi n/k)}
=-\sqrt{k}\left((-1)^{b/2}g(\chi)+E(b,\chi)\right).
\end{equation}
\end{corollary}

In particular, if $b=2$, then $E(2,\chi)=0$. We thus deduce the
following corollary of Corollary \ref{cor1}.

\begin{corollary}\label{cor2}
For $\chi$ as above,
\begin{equation}\label{e12}
\sum_{0<n<k/2}\chi(n)\df{\cos(2\pi n/k)}{\cos^2(\pi n/k)}
=g(\chi)\sqrt{k}.
\end{equation}
\end{corollary}

In particular,  from \eqref{e3}, we note that $g(\chi)$ is an
integer, and so  \eqref{e12} is an  analogue of Gauss's theorem
\eqref{gausssum} for even $\chi$, namely,
\begin{equation}\label{gausssumeven}
\sum_{n=1}^{k-1}\chi(n)\cos(2\pi n/k) =\sqrt{k}.
\end{equation}

Letting $k=5$ in \eqref{e12} and noting that $g(\chi)=4$, we find
that
\begin{equation}\label{e13}
\df{\cos(2\pi/5)}{\cos^2(\pi/5)}+\df{\cos(\pi/5)}{\cos^2(2\pi/5)}=4\sqrt{5}.
\end{equation}
Of course, since
\begin{equation}\label{cos}
\cos\left(\df{\pi}{5}\right)=\df{\sqrt{5}+1}{4}\qquad \text{and}
\qquad\cos\left(\df{2\pi}{5}\right)=\df{\sqrt{5}-1}{4},
\end{equation}
\eqref{e13} may  be easily verified directly.

\begin{corollary}\label{2cor1}
For $\chi$ as in Theorem \textup{\ref{thm1}}, we have
\begin{equation}\label{6e1}
S_2(k,\chi) := \sum_{0<n<k/2}\chi(n)\sec^2(\pi n/k) = -g(\chi)\sqrt{k}.
\end{equation}
\end{corollary}
\begin{proof}
Let $a=0$ in Theorem \ref{thm1}.  Then \eqref{6e1} follows,
 since $E(0,b,\chi)=0$ for all $b$.
\end{proof}
We conclude this section with a few evaluations.  Throughout the sequel, we
set 
$$ \chi_p(n) = \left(\df{n}{p}\right),$$
where the right-hand side above denotes the Legendre symbol. 
\begin{corollary}
We have
\begin{align*}
S_2(5,\chi) &= -4\sqrt{5},
&S_2(13,\chi) = -20\sqrt{13}, \\
S_2(17,\chi) &= 24\sqrt{17}, &S_2(29,\chi) = -60\sqrt{29}.
\end{align*}
\end{corollary}
\begin{proof}
These evaluations follow immediately from Corollary \ref{2cor1}, 
with $\chi(n)=\chi_p(n)$, which is even when $p$ is  congruent to
$1\pmod 4$. 
\end{proof}

\section{A General Theorem}\label{main}

%{\allowdisplaybreaks
%$$\leqalignno{\frac{f^3(-q)}{f^3(-q^5)}&=\frac{\psi(q)}{\psi(q^5)}\times
 %     \frac{\psi^2(q)-5q\psi^2(q^5)}{\psi^2(q)-q\psi^2(q^5)}, &(i) \cr

If $k$ is an odd positive integer, let $\chi_o(n)$ be a real, odd,
nonprincipal, primitive character of period $k$, and let
$\chi_e(n)$ be a real, even, nonprincipal, primitive character of
period $k$.  Define
\begin{align*}
H_o(z)&:=\frac{G(z;\chi_o)}{G(\chi_o)}\frac{1}{e^{2\pi iz}-1},\\
H_e(z)&:=\frac{G(z;\chi_e)}{G(\chi_e)}\frac{1}{e^{2\pi iz}-1}.
\end{align*}

\begin{theorem}\label{MainThm}
Let $k$, $\chi_o$, $\chi_e$, $H_o$, and $H_e$ be as above, and let
$f_o(z)$ and $f_e(z)$ be meromorphic functions satisfying the
following properties:
\begin{enumerate}
\item[(i)] $f_o(z)=f_o(z+k)$ and $f_e(z)=f_e(z+k)$ for all $z\in\C$,
\item[(ii)] $f_o(-z)=-f_o(z)$ and $f_e(-z)=f_e(z)$ for all $z\in\C$,
\item[(iii)] $f_o(z)$ and $f_e(z)$ are analytic for $\textup{Im }(z)\neq 0$ as well
 as at integers not divisible by $k$,
\item[(iv)] $f_o$ and $f_e$ have poles at the  points $z^o_1,...,z^o_M$
and $z^e_1,...,z^e_{M'}$, respectively,
\item[(v)] $f_o(z)$ and $f_e(z)$ have Fourier expansions of the form
\begin{align}\label{oddexp}
f_o(z) &= \sum_{m=-D}^\infty A_m e^{2\pi i
mz/k}, \\
\label{evenexp} f_e(z) &= \sum_{m=-D'}^\infty B_m e^{2\pi i mz/k}
\end{align} in the upper half-plane.
\end{enumerate}
 Then,
 %if $D<k$ and $D'<k$,
\begin{equation}
 \label{ThmOdd} \sum_{0<n<k/2} \chi_o(n)f_o(n) =
-i\sqrt{k}\sum_{m=0}^D \chi_o(m)A_{-m} - \pi i\sum_{m=1}^M
R_{z^o_m}(H_of_o) - \pi iR_0(H_of_o) \end{equation}
 and
\begin{equation} \label{ThmEven}
\sum_{0<n<k/2} \chi_e(n)f_e(n) = \sqrt{k}\sum_{m=0}^{D'}
\chi_e(m)B_{-m} - \pi i\sum_{m=1}^{M'}R_{z^e_m}(H_ef_e) - \pi
iR_0(H_ef_e).
\end{equation}
\end{theorem}

\begin{proof}
 We  prove \eqref{ThmOdd}; the proof of \eqref{ThmEven} is similar.

  Note that
\begin{align*}
H_o(-z) &= \frac{\sum_{j=1}^{k-1}\chi_o(j)e^{2\pi ij(-z)/k}}{G(\chi_o)}\frac{-e^{2\pi iz}}{e^{2\pi iz}-1} \\
&=\frac{-\sum_{j=1}^{k-1}\chi_o(j)e^{2\pi i(k-j)z/k}}{G(\chi_o)(e^{2\pi iz}-1)} \\
&=\frac{\sum_{j=1}^{k-1}\chi_o(k-j)e^{2\pi i(k-j)z/k}}{G(\chi_o)(e^{2\pi iz}-1)} \\
&= H_o(z).
\end{align*}
Thus, the function
\begin{equation*}
F_o(z) := H_o(z)f_o(z)
\end{equation*}
is an odd function with period $k$. This function has poles at
$z^o_1,\dots,z^o_M$, and (at most) simple poles at the points
$z=1,\dots,k-1$.  It also may have a pole at $z=0$. Let $C_N$
denote the same positively oriented indented rectangle
% with
%horizontal sides passing through $\pm iN$ and vertical sides
%passing through $0$ and $k$.
as in Theorem \ref{thm1}. Hence, by the residue theorem,
\begin{equation}
\label{MainInt} \frac{1}{2\pi i}\int_{C_N} F_o(z)dz =
\sum_{j=1}^{k-1}R_j(F_o) + \sum_{m=1}^{M}R_{z^o_m}(F_o) +
R_0(F_o). \ee

We first compute the residues at the integers $z=1,\dots,k-1$.
Since these are (at most) simple poles, the residue at $z=j$ is
given by

\begin{equation*}
R_j(F_o) = \lim_{z\ra j}(z-j)F_o(z) = \lim_{z\ra j}
\frac{G(z;\chi_o)}{G(\chi_o)}f_o(z)\frac{z-j}{e^{2\pi iz}-1} =
\frac{1}{2\pi i}\chi_o(j)f_o(j),
\end{equation*}
where we have used \eqref{factor} in the last step. By the oddness
and periodicity of $\chi_o$ and $f_o$, we find that
\begin{equation}\label{resintegers}
\sum_{j=1}^{k-1}R_j(F_o) =\sum_{j=1}^{k-1} \df{1}{2\pi
i}\chi_o(j)f_o(j) = \df{1}{\pi i}\sum_{1\leq j < k/2}
\chi_o(j)f_o(j).
\end{equation}

Secondly, we evaluate the integral on the left-hand side of
\eqref{MainInt} directly.  By periodicity, the integrals along the
vertical sides of $C_N$ cancel.   Thus, using the oddness of
$F_o(z)$ and periodicity, we see that the integral along the upper
horizontal edge is equal to
\begin{equation}
\int_{z=k+iN}^{z=iN} F_o(z) dz = -\int_{-z=-k-iN}^{-z=-iN}
F_o(-z)dz
%&= \int_{-z=-k-iN}^{-z=-iN} F_o(-z)d(-z) \notag\\
= \int_{u=-k-iN}^{u=-iN} F_o(u)du = \int_{u=-iN}^{u=k-iN}
F_o(u)du,\label{edge}
\end{equation}
which is the integral along the lower horizontal edge.  Also,
since the right-hand side of \eqref{MainInt} is independent of
$N$, we may let $N$ tend to $\infty$.  Thus, it remains  to
evaluate
\begin{equation*}
 \lim_{N\ra\infty}\frac{1}{2\pi
i}\int_{C_N}F_o(z)dz=\lim_{N\ra\infty}\frac{1}{\pi i}\int_k^0
F_o(x+iN) dx.
\end{equation*}

Now let \eb\label{mu} \mu:=e^{2\pi iz/k} \ee
 and expand $F_o(z)$
in a power series in  $\mu$.  Since
\begin{equation*}
\frac{1}{e^{2\pi iz}-1} = - 1 - \mu^k - \mu^{2k} - \cdots,
\end{equation*}
\begin{align}
F_o(z) &=
\frac{1}{G(\chi_o)}\left(\sum_{j=0}^{k-1}\chi_o(j)\mu^j\right)
 \left(\sum_{m=-D}^\infty A_m\mu^m\right)\left(-\sum_{n=0}^\infty \mu^{kn}\right) \nonumber \\
 &= \frac{1}{G(\chi_o)}\left(-\sum_{n=0}^\infty \sum_{j=0}^{k-1} \chi_o(j)\mu^{kn+j}\right)
 \left(\sum_{m=-D}^\infty A_m\mu^m\right) \nonumber \\
 &= \frac{1}{G(\chi_o)}\left(-\sum_{n=0}^\infty \chi_o(n)\mu^{n}\right)
 \left(\sum_{m=-D}^\infty A_m\mu^m\right) \nonumber \\
\label{intsum} &=\frac{1}{G(\chi_o)}\sum_{m=-D}^\infty c_m\mu^m
\end{align}
for some constants $c_m, m\geq -D$.  Since
\begin{equation*}
\int_k^0 \sum_{m=1}^\infty c_m\mu^m dx
\end{equation*}
 tends to 0
as $N$ tends to $\infty$, and
\begin{equation*}
\int_k^0 \mu^m dx=0
\end{equation*}
 for any $m\neq 0$, we find that, upon the use of \eqref{edge},
%and the inequality $D<k$,
\begin{equation}\label{intCN}
\lim_{N\ra\infty}\df{1}{2\pi i}\int_{C_N} F_o(z)dz =
-\frac{kc_0}{\pi iG(\chi_o)}=\frac{k}{\pi iG(\chi_o)}\sum_{m=0}^D
A_{-m}\chi_o(m),
\end{equation}
by \eqref{oddexp}.

In summary, we utilize \eqref{resintegers} and \eqref{intCN} in
\eqref{MainInt} to conclude that
\begin{equation}\label{MainFinal}
\frac{k}{\pi iG(\chi_o)}\sum_{m=0}^D A_{-m}\chi_o(m)=\df{1}{\pi i}
\sum_{1\leq j < k/2} \chi_o(j)f_o(j)+ \sum_{m=1}^{M}R_{z^o_m}(F_o)
+ R_0(F_o).
\end{equation}
Upon using the evaluation \eqref{gausssum} for $G(\chi_o)$, we see
that \eqref{MainFinal} is equivalent to \eqref{ThmOdd}.
\end{proof}

\section{Two Theorems on Trigonometric Sums and Corollaries}

We observe that for any integer $b$ the sine and cosine functions
satisfy the following properties:
\begin{align*}
\sin(b\pi(x+k)/k) &= (-1)^b\sin(b\pi x/k),\\
\cos(b\pi(x+k)/k) &= (-1)^b\cos(b\pi x/k).
\end{align*}

Thus, we may construct odd and even functions satisfying
properties (i)--(v) in Theorem \ref{MainThm} by taking appropriate
products and quotients of sines and cosines.  In particular, we
have the following theorem.

\begin{theorem}\label{EvenCor}
Let $k$ be an odd positive integer, and let $\chi_e$ be a real,
even, nonprincipal, primitive character of period $k$. Let $L$,
$a$, and $J$ be nonnegative integers with $a\leq J+1$.  Let
$b_1,\dots,b_L, c_1,\dots,c_L, d_1,\dots,d_J$ be positive integers
such that the sum
$$E:=-a+\sum_{\ell=1}^L(b_\ell -
c_\ell)+\sum_{j=1}^J d_j$$ is even, and that $d_1,\dots,d_{a-1}$
are odd.  Suppose also that $(c_\ell,k)=1$ for $1\leq\ell\leq L$,
and  $(c_i,c_j)=1$ for each $i\neq j$.
% Also, suppose that $E/2<k$.
Define \eb\label{g-def} g(\chi_e) = \sum_{j=1}^{k-1} (-1)^j
j\chi_e(j), \ee and let $P_o(n)$ \textup{(}resp.,
$P_e(n)$\textup{)} denote the number of
 solutions in the $(2L+J+a)$-tuples \\
$(\varepsilon_1,\dots,\varepsilon_L,m_1,\dots,m_L,\varepsilon_1',\dots,
\varepsilon_J',m_1',\dots,m_a')$ to the equation \eb \label{Reps}
n=\varepsilon_1b_1+\cdots+\varepsilon_Lb_L+m_1c_1+\cdots+m_Lc_L
+\varepsilon_1'd_1 +\cdots+\varepsilon_J'd_J+m_1'+\cdots+m_a', \ee
 where
$\sum\varepsilon_i+\sum m_i'$ is odd \textup{(}resp.,
even\textup{)}, and $\varepsilon_i, \varepsilon_i' \in \{0,1\},
m_i, m_i'\in\N\cup\{0\}$. Then {\allowdisplaybreaks\begin{align}
&\sum_{0<n<k/2} \chi_e(n)\left(\prod_{\ell=1}^L
\frac{\sin(b_\ell\pi n/k)}{\sin(c_\ell\pi n/k)}\right)
\left(\frac{\prod_{j=1}^J\cos(d_j\pi n/k)}{\cos^a(\pi n/k)}\right)
\nonumber
\\
& = \sqrt{k}\bigg\{2^{a-J}\sum_{m=0}^{E/2}
\chi_e(m)\left(P_e\left(\frac{E}{2}-m\right)-P_o\left(\frac{E}{2}-m\right)\right) \nonumber \\
& \quad- i\bigg(\sum_{M=1}^L\sum_{\substack{1\leq n < c_M \\
n\neq c_M/2}}
\frac{(-1)^n}{c_M}\frac{G(nk/c_M,\chi_e)\prod_{\ell=1}^L\sin(nb_\ell\pi/c_M)
\prod_{j=1}^J\cos(nd_j\pi/c_M)}{(e^{2\pi
ink/c_M}-1)\cos^a(n\pi/c_M)
\prod_{\ell\neq M}\sin(nc_\ell\pi/c_M)}\bigg) \nonumber \\
&\label{Cor1}  \quad- \sum_{\substack{1\leq M \leq L \\ c_M \text{
even}}} \frac{(-1)^{c_M/2}}{c_M}\frac{(-1)^a
g(\chi_e)\prod_{j=1}^{a-1}
(-1)^{(d_j-1)/2}d_j\prod_{\ell=1}^L\sin(b_\ell\pi/2)\prod_{j=a}^J
\cos(d_j\pi/2)}{\prod_{\ell\neq M}\sin(c_\ell\pi/2)}\bigg\}.
\end{align}}
\end{theorem}

In the preceding theorem, and throughout the rest of this section,
we adopt the convention that the empty products $\prod_{j=1}^0
a_j=1$ and $\prod_{j=1}^{-1} a_j=0$.

\begin{proof} In Theorem \ref{MainThm}, let $f_e(z)$ be defined by
\eb f_e(z) = \left(\prod_{\ell=1}^L \frac{\sin(b_\ell\pi
z/k)}{\sin(c_\ell\pi
z/k)}\right)\left(\frac{\prod_{j=1}^J\cos(d_j\pi z/k)}{\cos^a(\pi
z/k)}\right). \ee The function $f_e$ has possible poles at the
points $nk/c_\ell$, where $1\leq n<c_\ell$. We note that none of
these values are integers, because $(c_\ell,k)=1$ for all $\ell$.
Since $(c_i,c_j)=1$ for all $i\neq j$, these values are distinct.
Thus, these are at most simple poles, with the possible exception
of $k/2$, which occurs if one of the $c_\ell$ is even. Since at
least $a-1$ of the $d_j$ are odd, the point $z=k/2$ contributes a
pole of order at most 1 from the cosine factors. Since the numbers
of sine factors in the numerator and denominator are equal,  $z=0$
is a removable singularity.

Since $H_e$ has a removable singularity at $z=0$ and a simple zero
at $z=k/2$ (since $G(k/2,\chi_e)=0$), we conclude that $H_ef_e$
has at most simple poles at $nk/c_\ell$ for each $1\leq n< c_\ell,
1\leq\ell\leq L$.  Thus, the residue at the point $nk/c_M$, where
$1\leq M\leq L$ and $n\neq c_M/2$, is
\begin{align}
R_{nk/c_M}(H_ef_e) &=\lim_{z\ra nk/c_M}
\frac{z-nk/c_M}{\sin(c_M\pi z/k)}\frac{G(z,\chi_e)
\prod_{\ell=1}^L\sin(b_\ell\pi z/k)\prod_{j=1}^J\cos(d_j\pi
z/k)}{G(\chi_e)(e^{2\pi iz}-1)
\cos^a(\pi z/k)\prod_{\ell\neq M}\sin(c_\ell\pi z/k)} \nonumber \\
&= \frac{(-1)^n
k}{c_M\pi}\frac{G(nk/c_M,\chi_e)\prod_{\ell=1}^L\sin(nb_\ell\pi/c_M)
\prod_{j=1}^J\cos(nd_j\pi/c_M)}{G(\chi_e)(e^{2\pi
ink/c_M}-1)\cos^a(n\pi/c_M)\prod_{\ell\neq
M}\sin(nc_\ell\pi/c_M)}.\label{4.5}
\end{align}

When $z=(c_M/2)k/c_M$, where $c_M$ is the unique even $c_\ell$ (if
there is such a $c_\ell$), the residue is
\begin{align}
R_{k/2}(H_ef_e) &=\lim_{z\ra k/2} \frac{z-k/2}{\sin(c_M\pi
z/k)}\frac{G(z,\chi_e) \prod_{j=1}^{a-1}\cos(d_j\pi
z/k)}{\cos^a(\pi z/k)}\notag \\&\quad \times\frac{\prod_{\ell=1}^L
\sin(b_\ell\pi z/k)\prod_{j=a}^J\cos(d_j\pi
z/k)}{G(\chi_e)(e^{2\pi iz}-1)
\prod_{\ell\neq M}\sin(c_\ell\pi z/k)} \nonumber \\
&= \frac{(-1)^{c_M/2}k}{c_M\pi}\frac{2\pi
ig(\chi_e)\prod_{j=1}^{a-1}
(-1)^{(d_j-1)/2}d_j\prod_{\ell=1}^L\sin(b_\ell\pi/2)
\prod_{j=a}^J\cos(d_j\pi/2)}{\pi
(-1)^aG(\chi_e)(-2)\prod_{\ell\neq
M}\sin(c_\ell\pi/2)}.\label{4.6}
\end{align}

Lastly, we need to compute the Fourier expansion of $f_e(z)$ in
the upper half-plane.  We note that, with $\mu=e^{2\pi iz/k}$ as
in \eqref{mu},
\begin{align}
f_e(z) &= \left(\prod_{\ell=1}^L\frac{e^{ib_\ell\pi
z/k}-e^{-ib_\ell\pi z/k}} {e^{ic_\ell\pi z/k}-e^{-ic_\ell\pi
z/k}}\right)2^{a-J} \frac{\prod_{j=1}^J (e^{id_j\pi
z/k}+e^{-id_j\pi z/k})}{(e^{i\pi z/k}
+e^{-i\pi z/k})^a} \nonumber \\
&= 2^{a-J}\mu^{-\{\sum_{\ell=1}^L (b_\ell - c_\ell) - a +
\sum_{j=1}^J
d_j\}/2}\prod_{\ell=1}^L\frac{1-\mu^{b_\ell}}{1-\mu^{c_\ell}}
\frac{\prod_{j=1}^J(1+\mu^{d_j})}{(1+\mu)^a} \nonumber \\
&=
2^{a-J}\mu^{-E/2}\prod_{\ell=1}^L(1-\mu^{b_\ell})\prod_{\ell=1}^L
\left(\sum_{m=0}^\infty
\mu^{mc_\ell}\right)\prod_{j=1}^J(1+\mu^{d_j})
\left(\sum_{m=0}^\infty(-1)^m\mu^m\right)^a.\label{4.7}
\end{align}
By \eqref{ThmEven}, we need to determine the coefficient $B_{-m}$
of $\mu^{-m}$, as defined in \eqref{evenexp}. Collecting powers in
\eqref{4.7}, we see that
\begin{equation}\label{4.75}
B_{-m} = 2^{a-J}\left(P_e\left(\frac{E}{2}-m\right) -
P_o\left(\frac{E}{2}-m\right)\right).
\end{equation}
 Using \eqref{4.5}, \eqref{4.6}, and \eqref{4.75} in \eqref{ThmEven}, along with the evaluation of $G(\chi_e)$ in
\eqref{gausssum}, we complete the proof.
\end{proof}

\begin{theorem}\label{thm2}
For each pair of odd positive integers $a$ and $b$,
%with $ab-3<2k$,
 \eb \sum_{0<n<k/2} \chi_e(n)\frac{\cos^a(b\pi
n/k)}{\cos(\pi n/k)} = \sqrt{k}F(a,b,\chi),  \ee where
\begin{equation*}
F(a,b,\chi) :=\df{1}{2^{a-1}}\sum_{\substack{n,j,r \geq 0 \\
2(n+j+br)=ab-1}}(-1)^j\chi(n)\binom{a}{r},
\end{equation*}
where the sum is over all nonnegative integers $n$, $j$, and $r$
satisfying the condition $2(n+j+br)=ab-1$.
\end{theorem}

\begin{proof} In Theorem \ref{EvenCor}, let $L=0$, $a=1$, $J=a$, and
$d_1=\dots=d_J=b$.  Then we find that $E=ab-1$ is even, and that
the second and third sums on the right-hand side of \eqref{Cor1}
equal $0$.  We also find that the right-hand side of \eqref{Reps}
becomes $b(\varepsilon_1'+\cdots+\varepsilon_a')+m_1'$. Therefore,
\begin{equation*}
P_e\left(\df{ab-1}{2} - m\right) - P_o\left(\df{ab-1}{2} -
m\right) = \sum_{\substack{m_1',r\geq 0 \\ (ab-1)/2-m-m_1'-br=0}}
(-1)^r\binom{a}{r},
 \end{equation*}
where the sum is over all nonnegative integers $m_1', r$
satisfying the condition \linebreak $(ab-1)/2-m-m_1'-br=0$.
Rearranging the expression for the indices of summation, and
summing over $m$, we complete the proof.
\end{proof}

We now examine some special cases of Theorem \ref{thm2} when
$b=1$.

\begin{corollary}\label{cor3} If $a$ is odd and $\chi$ is even,
\begin{equation}\label{e28}
\sum_{0<n<k/2}\chi(n)\cos^{a-1}(\pi n/k) = \sqrt{k}F(a,1,\chi).
\end{equation}
\end{corollary}

\begin{proof} Set $b=1$ in Theorem \ref{thm2}.
\end{proof}

If $a=1$, then trivially $F(1,1,\chi)=0$, and so Corollary
\ref{cor3} reduces to
\begin{equation}\label{easy}
\sum_{0<n<k/2}\chi(n) = 0,
\end{equation}
which is easy to establish directly.

If $a=3$, observe that $F(3,1,\chi)=\tfrac{1}{4}$. Thus, from
Corollary \ref{cor3},
\begin{equation}\label{e29}
\sum_{0<n<k/2}\chi(n)\cos^{2}(\pi n/k) = \tfrac{1}{4}\sqrt{k}.
\end{equation}
The evaluation \eqref{e29} is also elementary, because if one
applies the double angle formula for $\cos(2\pi n/k)$ on the left
side of \eqref{e29} and uses both \eqref{gausssumeven} and
\eqref{easy}, \eqref{e29} easily follows. If $k=5$ in \eqref{e29},
we deduce that
$$ \cos^2(\pi/5)-\cos^2(2\pi/5) = \tfrac{1}{4}\sqrt{5},$$
which of course is an easy consequence of \eqref{cos}.

If $a=5$, Corollary \ref{cor3} reduces to
\begin{equation}\label{e30}
\sum_{0<n<k/2}\chi(n)\cos^{4}(\pi n/k) =
\tfrac{1}{16}\left(4+\chi(2)\right)\sqrt{k}.
\end{equation}
Using the double angle for cosine twice, along with
\eqref{gausssumeven}, \eqref{easy}, and \eqref{e29}, we can verify
\eqref{e30} directly.

\begin{corollary}\label{4cor1}
Let $a$, $d$, and $J$ be nonnegative integers such that $a\leq
J+1$ and $E:=dJ-a+12$ is even.
%and $E/2<k$.
Then
\begin{align}
& \sum_{0<n<k/2} \chi_e(n) \frac{\sin(3\pi n/k) \sin(5\pi
n/k)\sin(7\pi n/k)\cos^J(d\pi n/k)}
{\sin^3(\pi n/k)\cos^a(\pi n/k)} \nonumber \\
&= \frac{\sqrt{k}}{2^{J-a}}\sum_{m=0}^{E/2} \chi_e(m)
\sum_{\substack{E/2-m-3\varepsilon_1-5\varepsilon_2-7\varepsilon_3-dj\geq 0
\\ \varepsilon_i\in\{0,1\},\, 0\leq j\leq J}}
(-1)^{j+\varepsilon_1+\varepsilon_2+\varepsilon_3}\notag\\
&\quad \times
\binom{E/2-m-3\varepsilon_1-5\varepsilon_2-7\varepsilon_3
-dj+2+a}{2+a}.\label{4cor1e}
\end{align}
\end{corollary}

\begin{proof} Put $L=3$ and set  $c_\ell=1$, $1\leq \ell \leq 3$.
Thus, the second and third sums on the right-hand side of
\eqref{Cor1} equal $0$.  By moving the terms involving
$\varepsilon_i$ and $\varepsilon_i'$ to the left-hand side of
\eqref{Reps}, we find that the number of representations in
$(m_1,m_2,m_3,m_1',...,m_a')$ of
$$E/2-m-\sum\varepsilon_ib_i-\sum
\varepsilon_i'd_i=m_1+m_2+m_3+m_1'+...+m_a'$$
 is equal to
$$\binom{E/2-m-\sum\varepsilon_ib_i-\sum
\varepsilon_i'd_i+3+a-1}{3+a-1}.$$
 The desired result now follows.
\end{proof}

\begin{corollary}\label{4cor2}
We have
\begin{align}
S_3(k)&:=\sum_{0<n<k/2} \chi_e(n)\df{\sin(3\pi n/k)
\sin(5\pi n/k)\sin(7\pi n/k)}{\sin^3(\pi n/k)} \notag \\
&=\sqrt{k}\sum_{m=0}^6 \chi_e(m)\left\{\binom{8-m}{2} - 
\binom{5-m}{2} - \binom{3-m}{2}\right\}.\label{4cor2e}
\end{align}
\end{corollary}

\begin{proof}
Let $a=J=d=0$ in Corollary \ref{4cor1}.  Therefore, $E=12$ is even. 
 The index $j$ in the inner sum on the right-hand
side of \eqref{4cor1e} is always $0$, and therefore
 $E/2-m-3\varepsilon_1-5\varepsilon_2-7\varepsilon_3-dj$ is
non-negative only when $(\varepsilon_1,\varepsilon_2,\varepsilon_3)
=(0,0,0)$, $(1,0,0)$, and $(0,1,0)$, leading us to
the expression on the right-hand side of \eqref{4cor2e}.
\end{proof}

\begin{corollary}
We have
\begin{align*}
S_3(13)&=13\sqrt{13}, &S_3(17)=19\sqrt{17}, \\
S_3(29)&=3\sqrt{29}, &S_3(37)=13\sqrt{37}.
\end{align*}
\end{corollary}

\begin{proof}
Apply Corollary \ref{4cor2} with $\chi_e(n)=\chi_p(n)$
 and $p=13$, $17$, $29$, and
$37$.
\end{proof}

The next result is an analogue of Theorem 7.1 in \cite{bz1}.
 Technically, this theorem  follows from Theorem
\ref{EvenCor}, but it is perhaps easier to derive the result from 
Theorem \ref{MainThm}.

\begin{theorem}\label{thm3}
Let $\chi$ be even, and suppose that $b$ is an even positive
integer. Then
\begin{equation}\label{e31}
%S_3(b,\chi) :=
\sum_{0<n<k/2}\chi(n)\sin(b\pi n/k)\cot(\pi n/k) =
\sqrt{k}H(b,\chi),
\end{equation}
where
\begin{equation}\label{e32}
H(b,\chi) := \dfrac{1}{2}\chi\left(\df{b}{2}\right) +
\sum_{\substack{n,j\geq 1\\2n+2j=b}}\chi(n).
\end{equation}
\end{theorem}

\begin{proof}
We let $f_e(z) = \sin(b\pi z/k)\cot(\pi z/k)$ and $\chi_e = \chi$
in Theorem \ref{MainThm}.  Observe that $f(z)$ has no poles at
nonintegral points.  To compute its Fourier expansion, note that
\begin{align}
\label{e36} \sin(b\pi z/k) &= -\df{1}{2i}\mu^{-b/2}(1-\mu^{b}) \\
\intertext{and}
 \label{e37} \cot(\pi z/k)&=
 -i\left(1+2\sum_{j=1}^\i\mu^{j}\right).
 \end{align}
Thus,
 \begin{align}
 f_e(z)&=\df{1}{2}\mu^{-b/2}(1-\mu^b)+\mu^{-b/2}(1-\mu^b)\sum_{j=1}^\i\mu^j
 \notag \\
 \label{e38}
 &=\df{\mu^{-b/2}-\mu^{b/2}}{2}+\mu^{-b/2}\sum_{j=1}^{b-1} \mu^j.
 \end{align}
Therefore the coefficient $B_{-m}$ of $\mu^{-m}$, when $m\geq 0$,
is
\begin{equation}\label{e39}
B_{-m} = \begin{cases} 1/2, \qquad &\text{if $m=b/2$},\\
1, \qquad &\text{otherwise}. \end{cases}
\end{equation}
Using \eqref{e39} in \eqref{ThmEven}, we obtain  \eqref{e31}, and
the proof of Theorem \ref{thm3} is complete.
\end{proof}

\begin{corollary}\label{cor4} For even $\chi$,
\begin{equation}\label{e42}
\sum_{0<n<k/2}\chi(n)\cos^2(\pi n/k) = \tfrac{1}{4}\sqrt{k}.
\end{equation}
\end{corollary}

\begin{proof}  Set $b=2$ and
use the identity $\sin(2\theta) =2\sin\theta\,\cos\theta$ in
\eqref{e31}. Then note that $H(2,\chi)=\tf{1}{2}$.
\end{proof}

The identity \eqref{e42} is identical to \eqref{e29}.

We now  derive an analogue of Theorem \ref{EvenCor} for odd
characters.  Recall that $h(-k)$ denotes the class number for the
imaginary quadratic field $\Q(\sqrt{-k})$.

\begin{theorem}\label{OddCor}
Let $k$ be an odd positive integer, and let $\chi_o$ be a real,
odd, nonprincipal, primitive character of period $k$. Let $L$,
$a$, and $J$ be nonnegative integers with $a\leq J+1$.  Let
$b_1,\dots,b_L, c_1,\dots,c_L, d_1,\dots,d_J$ denote positive
integers such that the sum
$$E':=-1-a+\sum_{\ell=1}^L(b_\ell -
c_\ell)+\sum_{j=1}^J d_j$$
 is even and such that
$d_1,\dots,d_{a-1}$ are odd.  Suppose also that $(c_\ell,k)=1$ for
$1\leq\ell\leq L$ and that $(c_i,c_j)=1$ for  $i\neq j$.
% Also, suppose that $E/2<k$.
Let $P_o'(n)$ \textup{(}resp., $P_e'(n)$\textup{)} denote the
number of solutions in the $(2L+J+a+1)$-tuples
$$(\varepsilon_1,\dots,\varepsilon_L,m_0,m_1,\dots,m_L,\varepsilon_1',
\dots,\varepsilon_J',m_1',\dots,m_a')$$ to the equation
 \eb \label{RepsOdd}
n=\varepsilon_1b_1+\cdots+\varepsilon_Lb_L+m_0+m_1c_1+\cdots+m_Lc_L
+\varepsilon_1'd_1+\cdots+\varepsilon_J'd_J+m_1'+\cdots+m_a', \ee
where $\sum\varepsilon_i+\sum m_i'$ is odd \textup{(}resp.,
even\textup{)}, $\varepsilon_i, \varepsilon_i' \in \{0,1\}$, and $
m_i, m_i'\in\N\cup\{0\}$. Then {\allowdisplaybreaks\begin{align}
&\sum_{0<n<k/2} \chi_o(n)\left(\frac{1}{\sin(\pi
n/k)}\prod_{\ell=1}^L \frac{\sin(b_\ell\pi n/k)}{\sin(c_\ell\pi
n/k)}\right)
\left(\frac{\prod_{j=1}^J\cos(d_j\pi n/k)}{\cos^a(\pi n/k)}\right) \nonumber \\
&  = \sqrt{k}\bigg\{-2^{1+a-J}\sum_{m=0}^{E'/2} \chi_o(m)
\left(P_e'\left(\frac{E'}{2}-m\right)-P_o'\left(\frac{E'}{2}-m\right)\right) \nonumber \\
&\quad- \sum_{M=1}^L\sum_{1\leq n <
c_M}\frac{(-1)^n}{c_M}\frac{G(nk/c_M,\chi_o)
\prod_{\ell=1}^L\sin(nb_\ell\pi/c_M)\prod_{j=1}^J\cos(nd_j\pi/c_M)}{(e^{2\pi
ink/c_M}-1) \sin(n\pi/c_M)\cos^a(n\pi/c_M)\prod_{\ell\neq
M}\sin(nc_\ell\pi/c_M)}
 \nonumber \\
& \quad - \frac{G(k/2,\chi_o)\prod_{j=1}^{a-1}
(-1)^{(d_j-1)/2}d_j\prod_{\ell=1}^L
\sin(b_\ell\pi/2)\prod_{j=a}^J\cos(d_j\pi/2)}{2\prod_{\ell=1}^L
\sin(c_\ell\pi/2)} \nonumber \\
&\label{Cor2} \quad +
%\frac{\sum_{j=1}^{k-1}j\chi_o(j)}{k}\prod_{\ell=1}^L\frac{b_\ell}{c_\ell}
h(-k)\prod_{\ell=1}^L\frac{b_\ell}{c_\ell} \bigg\}.
\end{align}}
\end{theorem}

\begin{proof} The proof is analogous to that of Theorem
\ref{EvenCor}. We highlight the main differences.

We let $f_o(z)$ in Theorem \ref{MainThm} be replaced by
\begin{equation}
f_o(z)
= \frac{1}{\sin(\pi z/k)}\left(\prod_{\ell=1}^L
\frac{\sin(b_\ell\pi z/k)}{\sin(c_\ell\pi
z/k)}\right)\left(\frac{\prod_{j=1}^J\cos(d_j\pi z/k)}{\cos^a(\pi
z/k)}\right).
\end{equation}

We now need to determine the poles of the function $H_of_o$.  Note
that again we have poles at each of the points $nk/c_\ell$, $1\leq
n < c_\ell$, $1\leq\ell\leq L$.  There is also (at most) a simple
pole at  $k/2$.  Note that for odd characters $G(k/2,\chi_o) \neq
0$, but there is a possible pole at $k/2$ arising from  an even
$c_\ell$. Finally, there exists a pole at $z=0$ arising from the
extra sine factor in the denominator.

For $1\leq M \leq L$, $1\leq n < c_M$,
\begin{align}
& \qquad R_{nk/c_M}(H_of_o)\notag\\ &=\lim_{z\ra nk/c_M}
\frac{z-nk/c_M}{\sin(c_M\pi z/k)}\frac{G(z,\chi_o)\prod_{\ell=1}^L
\sin(b_\ell\pi z/k)\prod_{j=1}^J\cos(d_j\pi
z/k)}{G(\chi_o)(e^{2\pi iz}-1) \sin(\pi z/k)\cos^a(\pi
z/k)\prod_{\ell\neq M}\sin(c_\ell\pi z/k)}
\nonumber \\
&= \frac{(-1)^n
k}{c_M\pi}\frac{G(nk/c_M,\chi_o)\prod_{\ell=1}^L\sin(nb_\ell\pi/c_M)
\prod_{j=1}^J\cos(nd_j\pi/c_M)}{G(\chi_o)(e^{2\pi
ink/c_M}-1)\sin(n \pi/c_M)\cos^a(n\pi/c_M)\prod_{\ell\neq
M}\sin(nc_\ell\pi/c_M)}.
\end{align}

Next,
\begin{align}
&\qquad  R_{k/2}(H_of_o)\notag\\
&=\lim_{z\ra k/2} \frac{(z-k/2)\prod_{j=1}^{a-1} \cos(d_j\pi
z/k)}{\cos^a(\pi z/k)}\frac{G(z,\chi_o)
\prod_{\ell=1}^L\sin(b_\ell\pi z/k) \prod_{j=a}^J\cos(d_j\pi
z/k)}{G(\chi_o)(e^{2\pi iz}-1)\sin(\pi z/k)
\prod_{\ell=1}^L\sin(c_\ell\pi z/k)} \nonumber \\
&= \frac{-k\prod_{j=1}^{a-1}
(-1)^{(d_j-1)/2}d_j}{\pi}\frac{G(k/2,\chi_o)
\prod_{\ell=1}^L\sin(b_\ell\pi/2)
\prod_{j=a}^J\cos(d_j\pi/2)}{G(\chi_o)(-2)\prod_{\ell=1}^L\sin(c_\ell\pi/2)}.
\end{align}

Finally,
\begin{align}
R_0(H_of_o) &=\lim_{z\ra 0} \frac{z}{\sin(\pi
z/k)}\frac{G(z,\chi_o)}{G(\chi_o)(e^{2\pi
iz}-1)}\left(\prod_{\ell=1}^L
\frac{\sin(b_\ell\pi z/k)}{\sin(c_\ell\pi z/k)}\right)\frac{\prod_{j=1}^J\cos(d_j\pi z/k)}{\cos^a(\pi z/k)} \nonumber \\
&=\frac{k}{\pi}\frac{\sum_{j=1}^{k-1}j\chi_o(j)}{kG(\chi_o)}\prod_{\ell=1}^L\frac{b_\ell}{c_\ell}
=-\frac{kh(-k)}{\pi
G(\chi_o)}\prod_{\ell=1}^L\frac{b_\ell}{c_\ell},
\end{align}
where we applied \eqref{class} in the last equality.

Lastly, we determine the Fourier expansion of $f_o(z)$.  We
readily find that, again with $\mu=e^{2\pi iz/k}$ as in
\eqref{mu},
\begin{align}
f_o(z) &= \frac{2i}{e^{i\pi z/k}-e^{-i\pi
z/k}}\left(\prod_{\ell=1}^L \frac{e^{ib_\ell\pi
z/k}-e^{-ib_\ell\pi z/k}}{e^{ic_\ell\pi z/k}-e^{-ic_\ell\pi
z/k}}\right)2^{a-J}
\frac{\prod_{j=1}^J (e^{id_j\pi z/k}+e^{-id_j\pi z/k})}{(e^{i\pi z/k}+e^{-i\pi z/k})^a} \nonumber \\
&= 2^{1+a-J}i\mu^{-\{-1+\sum_{\ell=1}^L (b_\ell - c_\ell) - a +
\sum_{j=1}^J d_j\}/2}\frac{-1}{1-\mu}\prod_{\ell=1}^L
\frac{1-\mu^{b_\ell}}{1-\mu^{c_\ell}}
\frac{\prod_{j=1}^J(1+\mu^{d_j})}{(1+\mu)^a} \nonumber \\
&= -2^{1+a-J}i\mu^{-E'/2}\left(\sum_{m=1}^\infty
\mu^m\right)\prod_{\ell=1}^L(1-\mu^{b_\ell})\prod_{\ell=1}^L
\left(\sum_{m=0}^\infty \mu^{mc_\ell}\right)\notag\\&\quad \times
\prod_{j=1}^J(1+\mu^{d_j})\left(\sum_{m=0}^\infty(-1)^m\mu^m\right)^a.
\end{align}
Collecting powers, we find that the coefficient $A_{-m}$ of
$\mu^{-m}$, as defined in \eqref{oddexp}, is equal to
$-2^{1+a-J}i\left(P_e'(\frac{E'}{2}-m)-P_o'(\frac{E'}{2}-m)\right)$.

If we integrate $F_e(z)$ over the contour $C_N$, apply the residue
theorem, and calculate directly the integral of $F_e(z)$ over
$C_N$, using the calculation of the preceding paragraph for the
horizontal sides, we complete the proof in the same manner as we
did for Theorem \ref{EvenCor}.
\end{proof}

Many of our corollaries to Theorem \ref{OddCor} require 
explicit values of $h(-k)$ \cite[p.~425]{bs}.
For convenience, we record those values that we need for examples, namely, 
\begin{equation}\label{classeval}
h(-7)=h(-11)=h(-19)=1, \qquad h(-23)=3.
\end{equation}

\begin{theorem}\label{cor5}
If  $b$ and $d$ are odd positive integers,
\begin{align}
S_4(b,d,k)&:=4\sum_{0<n<k/2} \chi_o(n)
\df{\sin(b\pi n/k)\cos(d\pi n/k)}{\sin^2(2\pi n/k)} \notag \\
\label{o1} &=\sqrt{k}\left(4I(b,d,\chi_o) +
\df{(-1)^{(b+d)/2}d}{2}G(k/2,\chi_o) + bh(-k)\right),
% - \df{b}{\sqrt{k}}\sum_{j=1}^{k-1}j\chi_o(j),
\end{align}
where
\begin{align}
%I(b,d,\chi_o) &:= \sum_{m,r\geq 0}
%\chi_o(m)\bigg\{(-1)^r(r+1)\left(\binom{(b+d)/2-m-r-1}{1}
%+\binom{(b-d)/2-m-r-1}{1}\right) \notag \\
%\label{o2} &\qquad+(-1)^{r+1}(r+1)\binom{(d-b)/2-m-r-1}{1}\bigg\}.
I(b,d,\chi_o) &:= \sum_{m\geq 0}
\chi_o(m)\bigg\{\sum_{r=1}^{(b+d-2m)/2}(-1)^{r}r\left(\df{b+d}{2}-m-r\right) \notag \\
&\qquad +\sum_{r=1}^{(b-d-2m)/2}(-1)^{r}r\left(\df{b-d}{2}-m-r\right) \notag \\
\label{o2} &\qquad+\sum_{r=1}^{(d-b-2m)/2} (-1)^{r+1}r\left(\df{d-b}{2}-m-r\right)\bigg\}.
\end{align}
\end{theorem}

\begin{proof}
Let $L=J=1$, $b_1=b$, $d_1=d$, $c_1=1$, and $a=2$ in Theorem
\ref{OddCor}. Then, upon the use of the identity
$\sin(2\theta)=2\sin\theta\cos\theta$, we find that the left-hand
side of \eqref{Cor2} is equal to $S_4(b,d)$.  We also find that
$E'=-1-2+(b-1)+d=b+d-4$ is even. We compute the right-hand side of
\eqref{Cor2}.  With our choice of parameters, the right-hand side
of \eqref{RepsOdd} becomes \eb
b\varepsilon_1+m_0+m_1+d\varepsilon_1'+m_1'+m_2'. \ee
 We let
$r=m_1'+m_2'$ and observe that for each fixed $r\geq 0$, we obtain
exactly $r+1$ ordered pairs of nonnegative integers $(m_1',m_2')$.
For each fixed choice of $\varepsilon_1$, $\varepsilon_1'$, and $r$,
there are exactly $E'/2 - m - (b\varepsilon_1+d\varepsilon_1'+r)+1$
choices for $m_0$, and that choice of $m_0$ fixes $m_1$.
Therefore, for each fixed $m$,
\begin{align}
&P_e'\left(\tf{1}{2}E'-m\right)-P_o'\left(\tf{1}{2}E'-m\right) \notag \\
&= \sum_{r=0}^{E'/2-m+1} (-1)^r(r+1)\left(\df{b+d-4}{2}-m-r+1\right) \notag \\
&\quad +\sum_{r=0}^{E'/2-d-m+1} (-1)^r(r+1)\left(\df{b-d-4}{2}-m-r+1\right) \notag \\
\label{o3} &\quad+\sum_{r=0}^{E'/2-b-m+1} (-1)^{r+1}(r+1)
\left(\df{d-b-4}{2}-m-r+1\right).
\end{align}
Note that the second sum on the right-hand side of \eqref{Cor2} is
empty, and thus is equal to zero.  Using \eqref{o3} in
\eqref{Cor2}, shifting the index $r$ to $r-1$, and simplifying the
 remaining sums with our choice of
parameters, we conclude the proof.
\end{proof}
\begin{corollary}
We have
\begin{align*}
S_4(3,1,7) &= 4\sqrt{7}, &S_4(1,3,7)=4\sqrt{7}, \\
S_4(3,1,11)&= 0, &S_4(1,3,11)=-8\sqrt{11}.
\end{align*}
\end{corollary}
\begin{proof}
We apply Corollary \ref{cor5} with $\chi_o(n)=\chi_p(n)$,
  which is odd when $p$ is a
prime congruent to $3\pmod 4$.  Using the values of $h(-7)$ and $h(-11)$ 
from \eqref{classeval}, and the evaluations
\begin{equation}
G(7/2,\chi_7)=2, \qquad G(11/2,\chi_{11})=-6,
\end{equation}
we readily obtain the values of $S_4$ for $(b,d,k)=(3,1,7)$, $(1,3,7)$, 
$(3,1,11)$, and $(1,3,11)$.
\end{proof}

\begin{corollary}\label{cor6}
If $b$ is a positive odd integer,
\begin{align}
S_5(b,k):=2\sum_{0<n<k/2} \chi_o(n) \df{\sin(2b\pi n/k)}{\sin^2(2\pi n/k)}
&=4\sqrt{k}\sum_{m\geq 0}\chi_o(m)\sum_{r=1}^{b-m}(-1)^{r-1}r(b-m-r) \notag \\
\label{o4} &\qquad - \df{b\sqrt{k}}{2}G(k/2,\chi_o)
%-\df{b}{\sqrt{k}}\sum_{j=1}^{k-1}j\chi_o(j).
+b\sqrt{k}h(-k).
\end{align}
\end{corollary}

\begin{proof}
Let $b=d$ in Theorem \ref{cor5}.
\end{proof}

\begin{corollary}
We have
\begin{alignat*}{2}
S_5(1,7)&=0, &\qquad S_5(1,11)&=8\sqrt{11}, \\
S_5(1,19)&=4\sqrt{19},&\qquad  S_5(1,23)&=0, \\
S_5(3,7)&=-4\sqrt{7},& \qquad  S_5(3,11)&=8\sqrt{11}, \\
S_5(3,19)&=8\sqrt{11},&\qquad  S_5(3,23)&=-4\sqrt{23}.
\end{alignat*}
\end{corollary}

\begin{proof}
Apply Corollary \ref{cor6} with $\chi_o(n)=\chi_p(n)$ 
for the pairs $(b,k)=(1,7)$, $(1,11)$, $(1,19)$,
$(1,23)$, $(3,7)$, $(3,11)$, $(3,19)$, and $(3,23)$.  We obtain the desired 
results upon using the values of $h(-k)$ from
\eqref{classeval} and noting that
\begin{alignat*}{2}
G(7/2,\chi_7)&=2, &\qquad G(11/2,\chi_{11})
&=-6,  \\
G(19/2,\chi_{19})&=-6, &\qquad G(23/2,\chi_{23})
&=6.
\end{alignat*}
\end{proof}

\begin{theorem}\label{cor7}
If $a$, $b$, and $d$ are integers with $a\geq 0$, $b\geq 0$ and
even, and $d$ odd, then
\begin{align}
&S_6(a,b,d,k):=2\sum_{0<n<k/2}\chi_o(n)
\df{\cos^{a+b}(d\pi n/k)}{\sin(2\pi n/k)\cos^a(\pi n/k)} \notag \\
&=\begin{cases}\sqrt{k}\left(-2^{2-b}J(a,b,d,\chi_o)
%- \df{1}{k}\displaystyle{\sum_{j=1}^{k-1}j\chi_o(j)}\right),
+h(-k)\right),
&\text{if $b>0$,} \\
\sqrt{k}\bigg(-4J(a,0,d,\chi_o) -
\df{(-1)^{a(d-1)/2}}{2}d^aG(k/2,\chi_o)
%- \displaystyle{\df{1}{k}\sum_{j=1}^{k-1}j\chi_o(j)}\bigg),
+h(-k)\bigg), \quad & \text{if $b=0$,}\end{cases}
\end{align}
where
\begin{equation}
J(a,b,d,\chi_o) := \sum_{\substack{m,r,m_0,n\geq 0 \\
2dn+a+2+2m+2r+2m_0=d(a+b)}}
\chi_o(m)(-1)^r\binom{a+r}{a}\binom{a+b}{n}.
\end{equation}
\end{theorem}

\begin{proof}
In Theorem \ref{OddCor}, replace $a$ by $a+1$ and let $L=0$,
$J=a+b$, and $d_1=\dots=d_J=d.$ Then upon applying the identity
$\sin(2\theta)=2\sin\theta\cos\theta$, we find that the left-hand
side of \eqref{Cor2} is equal to $S_5(a,b)$. Note that
$E'=-1-(a+1)+d(a+b) = a(d-1)+db-2$ is even. The right-hand side of
\eqref{RepsOdd} is
\begin{equation}
m_0+d(\varepsilon_1'+\cdots+\varepsilon_{a+b}')+m_1'+\cdots+m_{a+1}'.
\end{equation}
Therefore, it is easy to see that if we let
$r=m_1'+\cdots+m_{a+1}'$, for each fixed $m\geq 0$,
\begin{align}
&P_e'\left(\df{ad+bd-a}{2}-1-m\right)-P_o'\left(\df{ad+bd-a}{2}-1-m\right)
 \notag \\
&\qquad= \sum_{r\geq 0} (-1)^r\binom{a+r}{a}
\sum_{m_0\geq 0}\binom{a+b}{((ad+bd-a)/2-1-m-r-m_0)/d} \notag \\
\label{o4a} &\qquad=\sum_{\substack{r,m_0,n\geq 0 \\
2dn+a+2+2m+2r+2m_0=d(a+b)}}
 (-1)^r\binom{a+r}{a}\binom{a+b}{n}.
\end{align}
Using \eqref{o4a} and evaluating the remaining terms on the
right-hand side of \eqref{Cor2}, we easily complete the proof.
\end{proof}

\begin{corollary}\label{cor8}
If $a\geq 1$  is even, and if $b\geq 0$ is odd, then
\begin{align}
S_7(a,b,k)&:=2\sum_{0<n<k/2} \chi_o(n)\df{\cos^a(b\pi n/k)}{\sin(2\pi n/k)} 
\notag \\
&= \sqrt{k}\bigg(h(-k) -
2^{2-a}\sum_{\substack{m,r,m_0,n\geq 0 \\
 2(bn+m+r+m_0) = ab-2}} \chi_o(m)(-1)^r\binom{a}{n} \bigg).\label{cor8-1}
%-\df{1}{k}\displaystyle{\sum_{j=1}^{k-1}j\chi_o(j)}\bigg).
\end{align}
\end{corollary}

\begin{proof}
Let $a=0$, $b=a,$ and $d=b$ in Theorem \ref{cor7}.
\end{proof}

We note that Corollary \ref{cor8} is analogous to Theorem 5.5 of
\cite{bz1}.
 The main difference is that
the parameter $b$ in Corollary \ref{cor8} above is odd, whereas
the
 corresponding parameter in Theorem
5.5 of \cite{bz1} is even.

\begin{corollary}\label{cor9}
For odd $k$,
\begin{equation}
\sum_{0<n<k/2}\chi_o(n)\cot(\pi n/k) = \sqrt{k}h(-k).
\end{equation}
\end{corollary}
\begin{proof}
Let $a=2$, $b=1$ in Corollary \ref{cor8}.
\end{proof}
Corollary \ref{cor9} is a corrected version of \cite[Eq.~$(6.3)$]{bz1}, 
which is
Lebesgue's class number formula \cite{lebesgue}.

\begin{corollary}
We have
\begin{alignat*}{2}
S_7(4,1,7) &= \df{3\sqrt{7}}{4}, &\qquad S_7(4,1,11)&=\df{3\sqrt{11}}{4}, \\
S_7(8,3,7) &= \df{19\sqrt{7}}{64}, &\qquad S_7(8,3,11)&=\df{49\sqrt{11}}{64}.
\end{alignat*}
\end{corollary}
\begin{proof}
Let $\chi_o(n)=\chi_p(n)$ in Corollary \ref{cor8}, and 
evaluate the right hand side of \eqref{cor8-1}
for $(a,b,k)=(4,1,7)$, $(4,1,11)$, $(8,3,7)$, and $(8,3,11)$, using the
 values of $h(-7)$ and $h(-11)$ in \eqref{classeval}.
\end{proof}

\section{Evaluations of Trigonometric Sums Not Involving
Characters}\label{s3}

\begin{theorem}\label{thmb1}
Suppose that $a$, $b$, and $ k$ are positive integers, where $b>1$
and  $k$ is odd.  Then
\begin{equation}\label{b1}
S_8(a,b,k) :=  \sum_{ 0<n<k/2}  \frac{ \sin^a (  2 \pi b n / k  ) }
{ \sin^a ( 2 \pi n / k  ) } = - \frac 1 2 \, b^a +
  \frac{1}{2}kS(a,b),
  \end{equation}
where
\begin{equation}\label{b2}
S(a,b) := 2\sum_{\substack{ m,n,r\geq
0\\2bn+2m+rk=ab-a}}(-1)^n\binom{a}{n}\binom{a-1+m}{m},
\end{equation}
where  in the case $r=0$, the terms are to be multiplied by
$\tfrac{1}{2}$.
\end{theorem}

Theorem \ref{thmb1} includes four identities found by Liu 
\cite{zhiguo}, namely, the special cases $(a,b,k) = (1,2,7),$ $
(1,3,7)$, $(7,2,7)$, and $ (7,3,7)$.

\begin{proof}
Let
\begin{equation}\label{b3}
  f(z) :=  \frac{ \sin^a( 2 \pi bz/ k ) }{ \sin^a( 2 \pi z/ k)  } 
\cot \pi z \ ,
\end{equation}
and integrate over the same contour $C_N$ described at the
beginning of the proof of Theorem \ref{thm1}.  Observe that $f(z)$
has simple poles at $z=0, 1, 2, \dots , k-1$.  A simple
calculation shows that
\begin{equation}\label{b4}
R_0 = \df{b^a}{\pi}.
\end{equation}
 For $0 < n < k$, we easily find that
\begin{equation}\label{b5}
  R_n = \frac{ \sin^a (  2 \pi b n / k ) }{\pi \sin^a (  2 \pi n /k ) }
   = R_{ k-n }  .
\end{equation}
Hence, by the residue theorem, \eqref{b4}, and \eqref{b5},
\begin{equation}\label{integral1}
  \int_{ C_N } f(z)dz = 2ib^a + 4i \sum_{ 0<n<k/2}
   \frac{ \sin^a (  2 \pi b n / k ) }{ \sin^a (  2 \pi n /k )}.
\end{equation}

Next we calculate the integral above directly. From the
periodicity of $f(z)$, we see that the integrals along the
vertical sides of $C_N$ cancel.  To compute the integrals on the
horizontal pieces, let $\mu = e^{ 2 \pi i z/ k }$. Then
\begin{align*}
  \sin^a (  2 \pi bz/ k )
  &= \left( \frac i 2 \right)^a \mu^{ -ab } \sum_{ n=0 }^{ a } (-1)^n \binom a n  \mu^{ 2bn }, \\
  \sin^{ -a } (  2 \pi z/ k )
  &= \left( \frac i 2 \right)^{ -a } \mu^{ a } \sum_{ m = 0 }^\i \binom{ a-1+m }{ m } \mu^{ 2m }, \\
  \cot \pi z
  &= i \, \frac{ \mu^{ k } + 1 }{ \mu^{ k } - 1 } = -i \left( 1 + 2 \sum_{ j = 1 }^\i \mu^{ kj } \right).
\end{align*}
Hence,
\begin{equation}\label{b6}
  f(z) = -i \, \mu^{ -ab+a } \left( \sum_{ n=0 }^{ a }(-1)^n
  \binom a n  \mu^{ 2bn } \right) \left( \sum_{ m = 0 }^{\i} \binom{ a-1+m }{ m } \mu^{ 2m } \right)
   \left( 1 + 2 \sum_{ r =1 }^{\i} \mu^{ kr } \right) .
\end{equation}
As in the proofs in the preceding section, we need to determine
the constant term, say $C_1$,  in the expansion \eqref{b6}.  With
some care, we see that
\begin{equation}\label{b7}
C_1 = -iS(a,b),
\end{equation}
where $S(a,b)$ is defined by \eqref{b2}.  Hence, letting $N$ tend
to $\i$, we deduce that
\begin{equation}\label{b8}
\lim_{N\to\i}\int_{C_{NT}}f(z)dz = -\int_0^kC_1 = ikS(a,b).
\end{equation}
By an analogous argument, we also find that
\begin{equation}\label{b9}
\lim_{N\to\i}\int_{C_{NB}}f(z)dz = ikS(a,b).
\end{equation}

Finally, combining \eqref{integral1}, \eqref{b8}, and \eqref{b9},
we conclude that
\begin{equation*}
2ib^a + 4i \sum_{ 0<n<k/2}
   \frac{ \sin^a (  2 \pi b n / k ) }{ \sin^a (  2 \pi n /k )} =
   2ikS(a,b),
   \end{equation*}
   which is equivalent to \eqref{b1}
   \end{proof}

By  a similar argument, we can also prove the following theorem.

\begin{theorem}\label{thmb2}
Suppose that $a$, $b$, and $k$ are positive integers, where $b>1$, odd,
and  $k$ odd.  Then
\begin{equation}
S_9(a,b,k):=  \sum_{ 0<n<k/2}
\frac{ \cos^a (  2 \pi b n / k) }{ \cos^a (2 \pi n / k)} = - \frac
1 2  +
  \frac{1}{2}kT(a,b),
  \end{equation}
where
\begin{equation}
T(a,b) := 2\sum_{\substack{ m,n,r\geq 0\\2bn+2m+rk=ab-a}}
(-1)^m\binom{a}{n}\binom{a-1+m}{m},
\end{equation}
where in the case $r=0$, the terms are to be multiplied by
$\tfrac{1}{2}$.
\end{theorem}

We conclude our paper with several evaluations.  Note that $S_8(7,2,7)$ is
found in Liu's paper \cite{zhiguo}.
\begin{corollary}
We have
\begin{alignat*}{2}
S_8(7,2,7)&=-57, &\qquad S_8(7,2,11)&=-64, \\
S_8(7,2,13)&=-64, &\qquad S_8(7,2,19)&=-64, \\
S_9(7,3,7)&=-1369, &\qquad S_9(7,3,11)&=-2162, \\
S_9(7,3,13)&=-2555, &\qquad S_9(7,3,19)&=-3734.
\end{alignat*}
\end{corollary}

\begin{proof}
These evaluations follow immediately upon the use of Theorem \ref{thmb1}
 for $S_8$ and Theorem \ref{thmb2} for $S_9$.
\end{proof}

\end{document}